\documentclass[11pt]{article}
\usepackage{amsmath}
\usepackage{amssymb}
\usepackage{amsthm}
\numberwithin{equation}{section}
\newtheorem{theorem}{Theorem}[section]

\newtheorem{lemma}[theorem]{Lemma}

\newtheorem{proposition}[theorem]{Proposition}
\newtheorem{remark}[theorem]{Remark}

\newcommand\Proof{\noindent{\bf Proof}\quad}
\newfont{\got}{eufm10 scaled \magstep1}

\newcommand{\ZZ}{{\mathbb Z}}
\newcommand\CC{\mathbb{C}}
\newcommand{\RR}{{\mathbb R}}
\newcommand{\FSM}{{\cal M}}
\newcommand{\FSS}{{\cal S}}
\newcommand\bLP{\\[\bigskipamount]}
\newcommand\bPP{\\[\bigskipamount]\indent}
\newcommand\qbinom[3]{\genfrac[]{0pt}0{#1}{#2}_{#3}}
\newcommand\qbinomtext[3]{\genfrac[]{0pt}1{#1}{#2}_{#3}}
\newcommand\half{\frac12}
\newcommand\thalf{\tfrac12}
\newcommand\ii{{\rm i}}
\newcommand\al{\alpha}

\newcommand\ep{\epsilon}
\newcommand\tha{\theta}
\newcommand\ka{\kappa}
\newcommand\la{\lambda}
\newcommand\iy{\infty}
\newcommand\LHS{left-hand side}

\begin{document}
\title{On a $q$-extension of Mehta's eigenvectors \\
of the finite Fourier transform for $q$ a root of unity}

\author{Mesuma K. Atakishiyeva, Natig M. Atakishiyev\\
and Tom H. Koornwinder}

\date{}

\maketitle

\begin{abstract}
It is shown that the continuous $q$-Hermite polynomials for $q$ a
root of unity have simple transformation properties with respect to
the classical Fourier transform. This result is then used to construct
$q$-extended eigenvectors of the finite Fourier transform in terms of
these polynomials.
\medskip

\noindent 2000MSC: 33D45, 42A38
\end{abstract}

\section{Introduction}

The {\em finite Fourier transform}
\cite[Ch.~7]{SteinShak} (also called {\em discrete Fourier transform}) is
defined as the Fourier transform associated with the {\em finite abelian
group} $\ZZ_N:=\ZZ/(N\ZZ)$ {\em of integers modulo} $N$ \cite{Serre,LN},
just as the classical integral Fourier transform is the Fourier transform
associated with $\RR$. In concrete terms, it is a linear transformation
$\Phi^{(N)}$ of the space of functions on $\ZZ$ with period $N$ defined
by
\begin{equation}
(\Phi^{(N)}f)(r):=\frac1{\sqrt N}\,\sum_{s=0}^{N-1}
\exp\left(\frac{2\pi \ii}{N} \, rs\right)\,f(s),\qquad r\in\ZZ.
\label{K12}
\end{equation}
Equivalently, if we identify the $N$-periodic function $f$ on $\ZZ$ with
the vector
\[
\bigl(f(0),f(1),\ldots,f(N-1)\bigr)
\]
in $\CC^N$, then $\Phi^{(N)}$
is a unitary operator on $\CC^N$ with matrix elements
\begin{equation}
\Phi_{rs}^{(N)}=\frac{1}{\sqrt{N}}\,\exp
\left(\frac{2\pi \ii}{N} \, rs\right)\, ,\quad\quad
0\leq r,s \leq N-1 \,.                             \label{1.1}
\end{equation}
We are interested in a suitable basis of $N$ eigenvectors
$f_n^{(N)}$
($n=0,1,...,N-1$) of $\Phi^{(N)}$,
which thus should satisfy the eigenvalue equation
\begin{equation}
\sum_{s=0}^{N-1}\,\Phi_{rs}^{(N)}\,f_n^{(N)}(s) =
\lambda_n\,f_n^{(N)}(r)                                   \label{1.2}
\end{equation}
for suitable eigenvalues $\lambda_n$. Since the fourth power of
$\Phi^{\left(N \right)}$ is the identity operator (or matrix), the
$\lambda_n$'s can only be equal to $\pm 1$ or $\pm \ii$.

The finite Fourier transform has deep roots in classical pure mathematics
and it is also extremely useful in applications. See \cite{AusTol},
\cite{AudTer} for mathematical and historical details and \cite{BerEv},
\cite{BEW} for the relation with Gauss sums.

Mehta studied in \cite{Mehta} the eigenvalue problem (\ref{1.2}) and
found analytically a set of eigenvectors $F_n^{(N)}$ of the finite
Fourier transform $\Phi^{(N)}$ of the form
\begin{equation}
F_n^{(N)}(r):=\sum_{k=-\iy }^{\iy }\,
e^{-\frac{\pi}{N}\left( kN+r\right)^2} \,H_n\,
\Big(\sqrt{\tfrac{2\pi}{N}}\,(kN+r)\Big),\quad r=0,1,\ldots,N-1,  \label{1.3}
\end{equation}
where $H_n\left( x\right)$ is the Hermite polynomial of degree $n$ in $x$.
These eigenvectors $F_n^{(N)}$ correspond to the
eigenvalues $\lambda_n= \ii^n$, that is,
\begin{equation}
\sum_{s=0}^{N-1}\Phi_{rs}^{\left( N\right)}\,F_n^{(N)}(s)
=\ii^{n}\,F_n^{(N)}(r)\,, \quad \quad r=0,\ldots,N-1. \label{1.4}
\end{equation}
They can be considered as a discrete analogue of the well-known continuous
case where the Hermite functions
$e^{-x^{2}/2}\,H_n(x)$ are constant multiples of their own Fourier transforms:
\begin{equation}
\frac{1}{\sqrt{2\pi }}\int_{-\iy }^{\iy}e^{\ii xy-x^2/2}\,H_n(x)\,dx
=\ii^{n}\,e^{-y^2/2}\,H_n(y). \label{1.5}
\end{equation}
It was conjectured by Mehta \cite{Mehta} that the $F_n^{(N)}$
($n=0,1,\ldots,N-1$
if $N$ is odd and $n=0,1,\ldots,N-2,N$ if $N$ is even) are linearly
independent,
but this problem has not been any further considered until now. As shown by
Ruzzi \cite{Ruzzi}, these systems are in general not orthogonal.

In later work it was shown that many $q$-extensions of the classical
orthogonal polynomials satisfy simple transformation properties under
the Fourier transform (see \cite{ANM1} and references therein). Thus
it was natural to repeat Metha's construction of eigenfunctions of
the finite Fourier transform in that context. This was done in
\cite{AGR}, \cite{ANM2} and \cite{ARW}. In particular, it was shown
there that the finite Fourier transform provides a link between
continuous $q$-Hermite and $q^{-1}$-Hermite polynomials of Rogers,
as well as between families of Rogers-Szeg\H{o} and Stieltjes-Wigert
polynomials. It turned out that the same form of connection exists
also between discrete $q$-Hermite polynomials of types I and II, see
\cite{MAR}.

Zhedanov \cite{Zhed} considered continuous $q$-Hermite polynomials
for $q$ a root of unity and he obtained a discrete orthogonality on
finitely many points with complex weights for a finite system of such
polynomials. In the present paper we consider the transformation
properties of these polynomials under the integral and finite Fourier
transforms. Analytically, this turns out to be a straightforward
extension of the earlier results for the continuous $q$-Hermite
polynomials with $0<q<1$. However, the resulting formulas, a little
different from the case that $q$ is real, are interesting enough to
be displayed.

Another feature of the present paper, compared with \cite{AGR}, \cite{ANM2}
and \cite{ARW}, is that we emphasize a more conceptual approach by using
Theorem \ref{K13} due to Dahlquist \cite{Dahl} and Matveev \cite{Matv}
and the (trivial) Lemma \ref{K42}, rather than repeating a technical
argument in each special situation.

The contents of the paper are as follows. In section 2 we recall
some properties of the continuous $q$-Hermite polynomials for
general complex $q$, and in particular for $q$ a root of unity.
In section 3 we consider the behaviour of these polynomials
times a Gaussian under the integral Fourier transform. This
result, together with the Dahlquist-Matveev Theorem \ref{K13},
then gives a construction of functions behaving nicely under
the finite Fourier transform. In section 5, using Lemma \ref{K42}
we obtain from these functions and their Fourier images eigenfunctions
of the finite Fourier transform. Finally, section 6 concludes
the paper with a brief discussion of some further research
directions of interest.

Throughout our exposition we employ standard notations of the theory
of special functions (see, for example, \cite{AAR} and \cite{GR}).
%
%
\section{$q$-Hermite polynomials for $q$ a root of unity}
The continuous $q$-Hermite polynomials of Rogers (see \cite{Rog},
\cite{Allaw}, \cite{AskIsm}), denoted by $H_n\left(x|\,q\right)$,
can be generated for any $q\in\CC$ by the three-term recurrence
relation
\begin{equation}
2x\,H_n\left(x|\,q\right)= H_{n+1}\left(x|\,q\right)\,+\,(1-q^n)\,
H_{n-1}\left(x|\,q\right)                               \label{2.1}
\end{equation}
with initial condition $H_0\left(x|\,q\right)=1$. Their explicit form
as a finite Fourier series in $\tha$ ($x=\cos\tha$) is
given by
\begin{equation}
H_{n}\left(\cos\tha|\,q \right) = \sum_{k=0}^{n}\,\qbinom nkq
e^{\ii(n-2k)\theta}\;,                             \label{2.2}
\end{equation}
where the symbol $\qbinomtext nkq$ stands for the {\em $q$-binomial
coefficient}
\begin{equation}
\qbinom nkq:=\,\frac{(q;q)_n}{(q;q)_k\,(q;q)_{n-k}}
=\qbinom n{n-k}q\,.                     \label{2.3}
\end{equation}
Here $(a;q)_n$ is the $q$-shifted factorial, see \cite[(1.2.15)]{GR}.
The right-hand sides of \eqref{2.3} and \eqref{2.2} are well-defined
for all $q\in\CC$ because the $q$-binomial coefficients are polynomials
in~$q$. 

Since $\sin\tha=\cos(\thalf\pi-\tha)$ we can rewrite \eqref{2.2} as
\begin{equation}
H_{n}\left(\sin\tha|\,q \right)=\ii^n\,\sum_{k=0}^{n}\,\qbinom nkq
\,(-1)^k\,e^{\ii(2k-n)\theta}\;.                        \label{K8}
\end{equation}

The polynomials $H_n(x|\,q)$ are orthogonal polynomials for $0<q<1$.
For $q>1$ they are orthogonal polynomials in $\ii x$\;:
the $q^{-1}$-Hermite polynomials (see \cite{Ask}) denoted by
\begin{equation}
h_n(x|\,q):=\ii^{-n}\,H_n(\ii x|\,q^{-1}).
\label{K3}
\end{equation}

Another case of orthogonality, but not with positive weights, was
considered for $q$ a root of unity, see \cite{Zhed}. For $M$ a
positive integer put
\begin{equation}
q_{j,M}:=\exp{(2\pi\ii\,j/M)},\qquad j\in\{1,2,\ldots,M-1\}.\label{K4}
\end{equation}
For such $q=q_{j,M}$ \eqref{2.1} and \eqref{2.2} remain valid.
In particular, for $q=q_{j,M}$ with $j$ and $M$ co-prime and for
$n=M$, the only non-vanishing terms in \eqref{2.2} occur for
$k=0$ and $M$. Hence, for $j$ and $M$ co-prime, we have
\begin{equation}
H_M(\cos\tha|\,q_{j,M})\,=\,2\,\cos{M\theta}=:2\,T_M(\cos\tha)\,\,,
\label{2.7}
\end{equation}
where $T_M(x)$ is a Chebyshev polynomial of the first kind (see, for
example, \cite[Remark 2.5.3]{AAR}).

As pointed out in \cite{Zhed}, the polynomials $H_n(x|\,q_{j,M})$
($n=0,1,\ldots,M-1$) satisfy a discrete orthogonality with possibly
complex weights on the $M$ zeros of $T_M(x)$ if $j$ and $M$ are
co-prime. It is for functions suitably defined in terms of these
polynomials that we will discuss their integral and finite Fourier
transforms.

Finally observe that, by induction with respect to $m$ and $n$, we
derive from \eqref{2.1} and \eqref{2.7} that, for $j$ and $M$ co-prime,
\begin{equation}
H_{mM+n}(x|\,q_{j,M})=(2T_M(x))^m\,H_n(x|\,q_{j,M}),
\qquad n=0,1,\ldots,M-1,\quad m=0,1,\ldots\;.
\label{2.8}
\end{equation}
\section{Integral Fourier transform}
There are $q$-extensions of the eigenfunction result \eqref{1.5} for
the integral Fourier transform. These interrelate certain $q$-polynomial
families (see \cite{ANM1} and references therein). For the conti\-nuous
$q$-Hermite polynomials we obtain:
\begin{lemma}
The Fourier transform of the functions
$e^{-x^2/2}\,H_n(\sin(\la x)|\,q)$ is given by:
\begin{multline}
\frac{1}{\sqrt{2\pi}}\,\int_{-\iy}^{\iy}\,e^{\ii xy -\,x^2/2}
\,H_n\big(\sin(\la x) |\,q\big)\,dx = \ii^n\,e^{-n^2\la^2/2}
\,e^{-y^2/2}\\\times\sum_{k=0}^n\,\qbinom nk{q^{-1}}\,
\left(q^{-1}e^{-2\la^2}\right)^{k(k-n)}\,(-1)^k\,e^{-(2k-n)\la y}
,\quad \la,q\in\CC,\;q\ne0,1.                          \label{K6}
\end{multline}
In particular, if $q=e^{-2\la^2}$ then
\begin{equation}
\frac{1}{\sqrt{2\pi}}\,\int_{-\iy}^{\iy}\,e^{\ii xy -\,x^2/2}
\,H_n\big(\sin(\la x) |\,q\big)\,dx = q^{n^2/4}\,e^{-\,y^2/2}
\,H_n\big(\sin(\ii\la y)|\,q^{-1}\big).           \label{K31}
\end{equation}
\end{lemma}
\Proof
Substitute \eqref{K8} on the \LHS, take termwise Fourier transforms
(which turns down to the Fourier transform of the Gaussian), and use
that
\[
\qbinom nk{q^{-1}}\,=\,q^{k(k-n)}\,\qbinom nkq.
\]
If $q=e^{-2\la^2}$
then we see by \eqref{K8} that \eqref{K6} simplifies to \eqref{K31}.
\qed
\bPP
We consider two cases of equation \eqref{K31}. First let $0<q<1$ and define
$\ka$ by
\begin{equation}
q=\exp{(-2{\ka}^2)},\quad 0 < \ka < \iy.    \label{K7}
\end{equation}
Then, by \eqref{K3}, we obtain (see \cite{AtNag1}, \cite{ANM1}):
\begin{proposition}
The Fourier transform of the functions
$e^{-x^2/2} H_n(\sin \ka x |\,q)$ is given by:
\begin{equation}
\frac{1}{\sqrt{2\pi}}\,\int_{-\iy}^{\iy}\,e^{\ii xy -\,x^2/2}
\,H_n\left(\sin \ka x |\,q\right)\,dx = \ii^{n}\,q^{n^2/4}\,
e^{-\,y^2/2}\,h_n\left(\sinh \ka y|\,q\right).  \label{3.1}
\end{equation}
\end{proposition}

Second, let $q:=q_{j,M}$ as in \eqref{K4} and put $\la:=\al_{j,M}$, where
\begin{equation}
\al_{j,M}:=\sqrt{\frac{\pi j}M}\,e^{-\,\pi\ii/4},\quad
{\rm hence}\quad e^{-2\al_{j,M}^2}=q_{j,M}\quad{\rm and}\quad
\ii\,\al_{j,M}=\overline{\al_{j,M}}.               \label{K11}
\end{equation}
\begin{proposition}
\label{K28}
The Fourier transform of the functions
$e^{-x^2/2}\,H_n\big(\sin(\al_{j,M}\,x)|\,q_{j,M}\big)$ is given by:
\begin{equation}
\frac{1}{\sqrt{2\pi}}\int_{-\iy}^{\iy}
e^{\,\ii xy - \,x^2/2}\,H_n\big(\sin(\al_{j,M}\,x)|\,q_{j,M}\big)\,dx
=q_{j,M}^{\,n^2/4}\,e^{-\,y^2/2}\,
H_n\big(\sin(\overline{\al_{j,M}}\,y)|\,q^{-1}_{j,M}\big). \label{3.5}
\end{equation}
\end{proposition}
The Fourier inversion formula of \eqref{3.5} is just the result of
taking complex conjugates on both sides of \eqref{3.5}.
We will mostly work with \eqref{3.5} for $n=0,1,\ldots,M-1$,
but this formula remains valid for all nonnegative integer values of $n$.
In particular, for $n=mM$ ($m=0,1,\ldots$), formula \eqref{3.5} takes by
substitution of \eqref{2.8} the form
\begin{multline}
\frac{1}{\sqrt{2\pi}}\int_{-\iy}^{\iy}
\cos^m\left(M\pi/2-\sqrt{\pi jM}\,e^{-\ii\pi/4}\,x\right)\,
e^{\,\ii xy - \,x^2/2}\,dx\\
= \ii^{m^2jM}\,e^{-y^2/2}\,
\cos^m\left(M\pi/2-\sqrt{\pi jM}\,e^{\ii\pi/4}\,y\right).
\label{3.7}
\end{multline}
Depending on the value of $M\pmod4$ this may be further simplified.
\begin{remark}
\label{K36}
\rm
Formulas \eqref{3.1} and \eqref{3.5} can be considered as $q$-analogues
of formula \eqref{1.5}, in the sense that \eqref{1.5} can be obtained
as the limit for $q\uparrow1$ of \eqref{3.1} and as the limit for
$M\to\iy$ of \eqref{3.5} with $j$ fixed. Indeed, from \cite[(5.26.1)]{KS}
we have
\begin{equation}
\lim_{q\to1} \Big(\sqrt{\thalf(1-q)}\,\Big)^{-n}
H_n\Big(x\sqrt{\thalf(1-q)}\,\Big|\,q\Big)=H_n(x).
\label{K30}
\end{equation}
Since $\ka\sim((1-q)/2)^\half$ as $q\uparrow 1$, it follows from
\eqref{K30} that
\begin{align*}
\lim_{q\uparrow1}\,\ka^{-n}\,H_n\big(\sin(\ka x)|\,q\big)&=H_n(x),\\
\lim_{q\uparrow1}\,(\ii\ka)^{-n}\,H_n\big(\sin(\ii\ka y)|\,q^{-1}\big)&=H_n(y).
\end{align*}
Hence, in view of \eqref{K3}, equation \eqref{3.1} with both sides multiplied
by $\ka^{-n}$ tends to \eqref{1.5} as $q\uparrow1$.

As for \eqref{3.5} with $j$ fixed, we have $\al_{j,M}\sim((1-q_{j,M})/2)^\half$
as $M\to\iy$. Hence it follows from \eqref{K30} that
\begin{align*}
\lim_{M\to\iy}\,\al_{j,M}^{-n}\,H_n\big(\sin(\al_{j,M}\,x)|\,q_{j,M}\big)
&=H_n(x),\\
\lim_{M\to\iy}\,\al_{j,M}^{-n}\,
H_n\big(\sin(\overline{\al_{j,M}}\,y)|\,q_{j,M}^{-1}\big)&=H_n(y).
\end{align*}
Then \eqref{3.5} with $j$ fixed
and with both sides multiplied by $\al_{j,M}^{-n}$,
tends to \eqref{1.5} as $q\uparrow1$.
\end{remark}

\section{Finite Fourier transform}
It was observed by Dahlquist \cite[Theorem 1]{Dahl} and Matveev
\cite[Theorem 8.1]{Matv} that Mehta's result \eqref{1.4} is a special
case of the following more general relationship between integral Fourier
transform and finite Fourier transform, which can be obtained as an
immediate consequence of the Poisson summation formula. This result
will hold for $f$ in a wide class of functions on $\RR$, but for
convenience we only formulate it for $f\in\FSS$, the space of
Schwartz functions on $\RR$ (see \cite[Ch.~5, \S1.3]{SteinShak}).
\begin{theorem}
\label{K13}
Define a linear map $\FSM^{(N)}$ from $\FSS$ to the space of
$N$-periodic functions on~$\ZZ$ by
\begin{equation}
(\FSM^{(N)}f)(r):=\sum_{k\in\ZZ}\,f\Big(\sqrt{\tfrac{2\pi}{N}}\,(kN+r)\Big)\,,
\quad f\in\FSS,\;r\in\ZZ.
\label{K16}
\end{equation}
If $f,g\in\FSS$ are related by the integral Fourier transform
\begin{equation}
g(y)\,=\,\frac{1}{\sqrt{2\pi}}\,\int_{-\iy}^{\iy}\,e^{\ii xy}
\,f(x)\,dx \,,                                            \label{4.1}
\end{equation}
and if $F:=\FSM^{(N)}f$, $G:=\FSM^{(N)}g$ then $F$ and $G$ are related by
the finite Fourier transform \eqref{K12}:
\begin{equation}
G(r)=\frac1{\sqrt N}\,\sum_{s=0}^{N-1} \exp\left(\frac{2\pi \ii}{N}
\, rs\right)\,F(s).                                     \label{K14}
\end{equation}
\end{theorem}
In particular, in the case when $g(x)=\lambda f(x)$, $\lambda = \pm 1,
\pm \ii$, one has $G(r)=\lambda F(r)$; so Mehta's eigenvectors (\ref{1.3})
are a particular case of the general statement.

Let us now apply Theorem \ref{K13} to the case where $f$ and $g$
are implied by \eqref{K31}, i.e.,
\[
f(x):=e^{-x^2/2}\,H_n\big(\sin(\la x)|\,q\big),\qquad
g(y):=e^{-y^2/2}\,H_n\big(\sin(\ii\la y)|\,q^{-1}\big),
\]
where $\la,q\in\CC$, $q\ne0,1$, and $q=e^{-2\la^2}$. Put
\begin{align}
f_n^{(N)}(r|\,q)&:=(\FSM^{(N)}f)(r)=
\sum_{k=-\iy}^\iy\,e^{-\frac\pi N\,(kN+r)^2}
H_{n}\Big(\sin\Big(\la\sqrt{\tfrac{2\pi}N}\,(kN+r)\Big)
\,\Big|\,q\Big)\,,                        \label{K37}\\
g_n^{(N)}(r|\,q)&:=q^{-n^2/4}\,(\FSM^{(N)}g)(r)\nonumber \\
&\;= \sum_{k=-\iy}^\iy\,e^{-\frac\pi N\,(kN+r)^2}
H_{n}\Big(\sin\Big(\ii\la\sqrt{\tfrac{2\pi}N}\,(kN+r)\Big)
\,\Big|\,q^{-1}\Big)\,.                        \label{K38}
\end{align}
Then we obtain by Theorem \ref{K13}:
\begin{proposition}
\label{K32}
The finite Fourier transform of the functions $f_n^{(N)}$ is given by:
\begin{equation}
\Phi^{(N)}\big(f_n^{(N)}(\,.\,|\,q)\big)(r)=\sum_{s=0}^{N-1}
\exp\left(\frac{2\pi \ii}{N} \, rs\right)\,f_n^{(N)}(s|\,q)=
q^{n^2/4}\,g_n^{(N)}(r|\,q).
\label{K33}
\end{equation}
\end{proposition}

Just as we did in \S3 for \eqref{K31}, we consider two special cases of
\eqref{K33}. First, let $0<q<1$ and take $\la=\ka$ as in \eqref{K7}.
Then \eqref{K33} holds (see \cite{AGR} and \cite{ARW}) with
\begin{equation}
g_n^{(N)}(r|\,q)=
\ii^n\,\sum_{k=-\iy}^\iy\,e^{-\frac\pi N\,(kN+r)^2}
h_{n}\Big(\sinh\Big(\ka\,\sqrt{\tfrac{2\pi}N}\,(kN+r)\Big)\,
\Big|\,q\Big)\,.
\label{K34}
\end{equation}

Second, consider \eqref{K33} with $q=q_{j,M}$ and $\la=\al_{j,M}$ as
in \eqref{K11}. Then \eqref{K33} holds with
\begin{align}
f_n^{(N)}(r|\,q_{j,M})&=
\sum_{k=-\iy}^\iy\,e^{-\frac\pi N\,(kN+r)^2}
H_{n}\Big(\sin\Big(\al_{j,M}\,\sqrt{\tfrac{2\pi}N}\,(kN+r)\Big)\,
\Big|\,q_{j,M}\Big)\,,
\label{K39}\\
g_n^{(N)}(r|\,q_{j,M})&=
\sum_{k=-\iy}^\iy\,e^{-\frac\pi N\,(kN+r)^2} H_{n}\Big(\sin\Big
(\overline{\al_{j,M}}\,\sqrt{\tfrac{2\pi}N}\,(kN+r)\Big)\,
\Big|\,q_{j,M}^{-1}\Big)\, \nonumber \\ &=\overline{f_n^{(N)}(r|\,q_{j,M})}.
\label{K35}
\end{align}

By similar reasoning as in Remark \ref{K36} we see that \eqref{1.3}
is a limit of \eqref{K33} in these two cases (as $q\uparrow1$ in the
first case and as $M\to\iy$ in the second case). These limits are
formal because we have to take termwise limits for $f_n^{(N)}(s|\,q)$
and $g_n^{(N)}(r|\,q)$. Just as for the functions $F_n^{(N)}$ given
by \eqref{1.3}, we do not know if the functions $f_n^{(N)}$ defined
by \eqref{K37} are linearly independent.

Let $\ep=\pm1$. From \eqref{K16} we see that, if $f(-r)=\ep f(r)$ for
all $r$, then also $(\FSM^{(N)}f)(-r)=\ep (\FSM^{(N)}f)(r)$ for all $r$.
Also note from \eqref{2.1} that $H_n\left(-x|\,q\right)=(-1)^n H_n\left
(x|\,q\right)$. Hence, by \eqref{K37} and \eqref{K38},
\begin{equation*}
f_n^{(N)}(-r|\,q)=(-1)^n f_n^{(N)}(r|\,q),
\qquad
g_n^{(N)}(-r|\,q)=(-1)^n g_n^{(N)}(r|\,q).
\end{equation*}
Thus consider \eqref{K33} with $f_n^{(N)}$ and $g_n^{(N)}$ as in
\eqref{K39}, \eqref{K35}. Then take complex conjugates on both sides and
use the simple facts just mentioned above.
Then we obtain:
\begin{align}
\Phi^{(N)}\big(f_n^{(N)}(\,.\,|\,q_{j,M})\big)(r)&=
q_{j,M}^{n^2/4}\,\overline{f_n^{(N)}(r|\,q_{j,M})},
\label{K40}\\
\Phi^{(N)}\big(\overline{f_n^{(N)}(\,.\,|\,q_{j,M})}\big)(r)&=
(-1)^n\,q_{j,M}^{-n^2/4}\,f_n^{(N)}(r|\,q_{j,M}).
\label{K41}
\end{align}
\begin{remark}\rm
It is tempting to consider the case $q=q_{j,M}$ and $\la=\al_{j,M}$ of
\eqref{K33} with $M=N$, even if one loses then the possibility to take
the limit to \eqref{1.3} for $M\to\iy$. One might hope to arrive at some
transform acting on the polynomials $H_n(x|\,q_{j,N})$ ($n=0,1,\ldots,N-1$;
$j$ and $N$ co-prime) but only involving their values at the zeros of $T_N(x)$,
i.e., at the points involved in Zhedanov's \cite{Zhed} orthogonality
relations for these polynomials. At the moment we have no idea how to
proceed here.
\end{remark}
\section{Finite Fourier $q$-extended eigenvectors}
We now show that relations \eqref{K40}, \eqref{K41} enable us to construct
$q$-extensions of Mehta's eigenvectors \eqref{1.3} of the finite Fourier
transform when the deformation parameter $q$ is a root of unity.
For this we need the following trivial observation.
\begin{lemma}
\label{K42}
Let $V$ be a complex linear space. Let $f,g\in V$ and $a,b\in\CC$ such
that $\Phi f=a^2g$ and $\Phi g=b^2 f$. Then
\[
\Phi(bf\pm ag)=\pm ab\,(bf\pm ag).
\]
\end{lemma}

Combination of this lemma with \eqref{K40} and \eqref{K41} shows that
the functions
\[
\ii^n\,q_{j,N}^{-n^2/8}\,f_n^{(N)}(\,.\,|\,q_{j,M})\pm
q_{j,M}^{n^2/8}\,\,\overline{f_n^{(N)}(\,.\,|\,q_{j,M})}
\]
are eigenfunctions of $\Phi^{(N)}$ with eigenvalues $\pm\,\ii^n$.
So put
\begin{align*}
F_n^{(N)}(r|\,q_{j,M})&:=\Re\left(e^{\ii\pi n/8}\,q_{j,M}^{-n^2/8}\,
f_n^{(N)}(r|\,q_{j,M})\right),\\
G_n^{(N)}(r|\,q_{j,M})&:=\Im\left(e^{\ii\pi n/8}\,q_{j,M}^{-n^2/8}\,
f_n^{(N)}(r|\,q_{j,M})\right).
\end{align*}
Then
\begin{align}
\Phi^{(N)}\big(F_n^{(N)}(\,.\,|\,q_{j,M})\big)(r)&=
\ii^n\,F_n^{(N)}(r|\,q_{j,M}),\\
\Phi^{(N)}\big(G_n^{(N)}(\,.\,|\,q_{j,M})\big)(r)&=
-\,\ii^n\,G^{(N)}(r|\,q_{j,M}).
\end{align}

\section{Concluding comments and outlook}

We have demonstrated that the continuous $q$-Hermite polynomials for $q$
a root of unity are interrelated by the classical Fourier transform, see
\eqref{3.5}. Then the technique of constructing eigenvectors of
the finite Fourier transform, pioneered by Mehta \cite{Mehta}
and formulated in a more
systematic way by Dahlquist \cite{Dahl} and Matveev \cite{Matv}, has
been employed in order to construct $q$-extended eigenvectors of the
finite Fourier transform.

Quite naturally, it would be of considerable interest to find out whether
there are other families of $q$-polynomials for $q$ a root of unity, which
also possess such simple transformations properties with respect to the Fourier
transform and, consequently, give rise to other $q$-variants of Metah's 
eigenfunctions \eqref{1.3} of the finite Fourier transform.
The point is that the continuous $q$-Hermite
polynomials, considered in the present paper, belong to the lowest level in
the Askey hierarchy of basic hypergeometric polynomials \cite{KS}. Therefore
it will be natural to attempt to apply the same technique to the study of
other $q$-families, which depend on some additional parameters (and therefore
occupy the higher levels in the Askey $q$-scheme).

Finally, another direction for further study is connected with $q$-extensions
of the harmonic oscillator in quantum mechanics \cite{Arik}, \cite{Macf},
\cite{Bied}. We remind the reader that for proving the fundamental formula
(\ref{1.4}) for Mehta's eigenvectors $F^{\left(n\right)}$ of the finite
Fourier transform operator (\ref{1.1}) it is vital to use the simple
transformation property (\ref{1.5}) of the Hermite functions $H_n(x)
\,\exp(-x^2/2)$ with respect to the Fourier transform. Moreover, these
eigenvectors $F^{\left(n\right)}$ are actually built in terms of these
Hermite functions taken at the infinite set of discrete points
$x_r^{(k)}:=\sqrt{\frac{2\pi}{N}}\,\left(kN+r \right),\,\,0\leq r\leq N-1\,,
\,\,k\in{\mathbb Z}$ ({\it cf.}~\eqref{K16}). In other words, Mehta's
technique of proving (\ref{1.4}) is based on introducing a discrete analogue
of the quantum harmonic oscillator. It seems that $q$-extended eigenvectors
of the finite Fourier transform, constructed in the foregoing sections, can
be similarly viewed as discrete analogues of the $q$-harmonic oscillator of
Macfarlane and Biedenharn for $q$ a root of unity.

We plan to continue our studies in both of these directions. But a focus of
our attention will be on constructing an explicit form of a finite-difference
equation for suitable eigenvectors \eqref{1.2} of the finite Fourier transform
operator $\Phi^{(N)}$. From
the group-theoretic point of view this amounts to finding an adequate
finite representation, associated with the Heisenberg--Weyl group
\cite{Serre,LN}. Purely analytically this is reduced to the construction
of an explicit finite-difference operator $L$ which commutes with
$\Phi^{(N)}$, so that the eigenspaces of
$\Phi^{(N)}$ can be decomposed by spectral decomposition of $L$.

\medskip
\noindent
{\bf Acknowledgements}

\medskip

NMA would like to thank the Korteweg-de Vries Institute for Mathematics,
University of Amsterdam, for the hospitality extended to him during his
visit in April, 2008, when part of this work was carried out. The
participation of MKA in this work has been supported by the CONACYT
project No.25564, whereas NMA has been supported by the UNAM-DGAPA
project IN105008 "\'Optica Matem\'atica". The authors thank a referee
for helpful comments on an earlier version of this paper.

\quad
\begin{footnotesize}
\begin{quote}
Mesuma K. Atakishiyeva \\
Facultad de Ciencias, Universidad Aut\'onoma del Estado de Morelos \\
C.P.~62250 Cuernavaca, Morelos, M\'exico\\
email: {\tt mesuma@servm.fc.uaem.mx}
\bLP
Natig M. Atakishiyev \\
Instituto de Matem\'{a}ticas, Unidad Cuernavaca \\
Universidad Nacional Aut\'onoma de M\'exico\\
A.P.~273-3 Admon.3, Cuernavaca, Morelos, 62251 M\'{e}xico
\\
email: {\tt natig@matcuer.unam.mx}
\bLP
Tom H. Koornwinder \\
Korteweg-de Vries Institute, University of Amsterdam\\
P.O.\ Box 94248,
1090 GE Amsterdam, The Netherlands\\
email: {\tt T.H.Koornwinder@uva.nl}
\end{quote}
\end{footnotesize}

\end{document}